\documentclass[a4paper,10pt]{article}

\usepackage{amssymb}
\usepackage{amsfonts}
\usepackage{amsmath}
\usepackage{amsthm}
\usepackage{stmaryrd}
\usepackage{enumerate}
\usepackage{dsfont}
\usepackage{a4wide}
\usepackage{graphicx}
\usepackage{color}
\usepackage{multicol}

 

\DeclareMathOperator{\var}{Var}

\def\defby{=:}
\def\bydef{:=}
\def\build#1_#2^#3{\mathrel{\mathop{\kern 0pt#1}\limits_{#2}^{#3}}}
\def\tend#1#2#3{\build\hbox to 12mm{\rightarrowfill}_{#1\rightarrow#2}^{#3}}

\def\reff#1{(\ref{#1})}

\def\sig#1{\sigma(#1)}
\def\siga#1{\sigma(a_{#1})}

\def\B#1{B_{#1}}
\def\Ba#1{B_{a_{#1}}}
\def\Bbar#1{\overline{B}_{{#1}}}
\def\Babar#1{\overline{B}_{a_{#1}}}
\def\Ca#1{C_{a_{#1}}}
\def\p#1{p(#1)}
\def\pa#1{p(a_{#1})}
\def\gp#1#2{p^{#1}(#2)}
\def\ella{\ell(a)}
\def\alf#1{\alpha(#1)}
\def\l#1{\lambda_#1}
\def\m#1{\mu_#1}

\def\mua#1{\mu^{(a)}_#1}
\def\la#1{\lambda_{a_{#1}}}
\def\ma#1{\mu_{a_{#1}}}
\def\D#1{D_{#1}}
\def\Da#1{D_{a_{#1}}}
\def\De#1{\Delta^2_{#1}}
\def\Dea#1{\Delta^2_{a_{#1}}}

\def\a#1{\alpha_{#1}}
\def\b#1{\beta_{#1}}
\def\at#1{\widetilde{\alpha}_{#1}}
\def\bt#1{\widetilde{\beta}_{#1}}

\newcommand{\X}[1]{X({#1})} 
\newcommand{\x}[2]{X^{#1}({#2})}

\newcommand{\Xa}[1]{X^{(a)}({#1})}
\newcommand{\T}[2]{T_{#1\to#2}}
\newcommand{\Ta}[2]{T^{(a)}_{#1\to#2}}

\newcommand{\Tat}[2]{\widetilde{T}^{(a)}_{#1\to#2}}
\newcommand{\TR}[1]{T_{#1}}
\newcommand{\TRa}[1]{T^{(a)}_{#1}}

\newcommand{\E}[1]{\mathbb{E}\left[{#1}\right]}
\newcommand{\EE}[1]{\mathbb{E}\bigl[{#1}\bigr]}
\newcommand{\EEE}[1]{\mathbb{E}\Bigl[{#1}\Bigr]}
\newcommand{\pr}[1]{\mathbb{P}\left({#1}\right)}
\newcommand{\prr}[1]{\mathbb{P}\bigl({#1}\bigr)}

\newcommand{\pii}[1]{\pi\left({#1}\right)}
\newcommand{\piii}[1]{\pi\bigl({#1}\bigr)}

\newcommand{\piia}[1]{\pi^{(a)}\left({#1}\right)}

\newcommand{\qa}{Q_a}
\newcommand{\ra}{R_a}
\newcommand{\Ga}{\Gamma_a}
\newcommand{\Q}[1]{Q({#1})}

\newcommand{\linf}[1]{\xrightarrow[#1\to\infty]{}}
\newcommand{\ze}{\mbox{\tiny0}}

\title{Cut-off and Escape Behaviors\\ for Birth and Death Chains on Trees}
\author{O. Bertoncini}

\begin{document}
\newtheorem{theo}{Theorem}[section]
\newtheorem{prop}[theo]{Proposition}
\newtheorem{lem}[theo]{Lemma}
\newtheorem{coro}[theo]{Corollary}
\theoremstyle{definition}
\newtheorem{definition}[theo]{Definition}
\newtheorem{rem}[theo]{Remark}

\maketitle


\begin{abstract}
We consider families of discrete time birth and death chains on trees, and show that in presence of a drift toward the root of the tree, 
the chains exhibit cut-off behavior along the drift and escape behavior in the opposite direction.
\end{abstract}



\section{Introduction}

Although they originally come from different research fields and seem apparently to be very different phenomena, cut-off and escape behaviors 
have been related at the level of hitting times for birth and death chains on the line (see \cite{barberfer09,olivphd,berbarfer08}). 
Cut-off behavior refers to the famous cut-off phenomenon first discovered and studied by Aldous and Diaconis in the 80's 
(see \cite{AldousRW,DiaAlShCST,DiaAlSUT}), which is characterized by an almost deterministic asymptotical abrupt convergence to equilibrium for families of 
Markov processes. In opposition, escape behavior is usually associated to the exit from metastability for a system trapped in local minimum of the 
energy profile. In that case, the transition to equilibrium occurs at ``unpredictable'' exponential times (see \cite{CGOV}). \\

We refer to the introduction in \cite{barberfer09} for a comparative historical discussion of the phenomena. 
In this reference ---dealing with birth-and-death chains on the line with drift toward the origin--- cut-off and escape phenomena are characterized 
by their distinct hitting-time behavior. It is shown that, under suitable drift hypotheses, the chains exhibit both cut-off behavior toward zero and 
escape behavior for excursions in the opposite direction. Furthermore, as the evolutions are reversible, the law of the final escape trajectory coincides 
with the time reverse of the law of cut-off paths.\\

In the present work, we extend this study to birth-and-death chains on trees with drift toward the root of the tree. 
Our chains are the discrete time counterpart of the birth-and-death processes studied by Mart\'inez and Ycart in~\cite{YcartMart}. 
Under a ``uniform'' drift condition, they prove that hitting times of zero starting from any state $a$ (denoted by $\T{a}{\ze}$) exhibit cut-off 
at mean times when $a$ goes to infinity (Theorem 5.1 in~\cite{YcartMart}). Their drift condition is expressed as an exponential decaying tail of the 
invariant probability measure $\pi$, in the sense that, for all $a$, the quantity $\pii{\B a}/\pii{a}$ is uniformly bounded by some $K$. 
In comparison, our drift condition concerns only the branch which contains the state $a$ (see formulas \reff{def:Ka} and \reff{eq:Ka}), 
and we allow the upper bound $K_a$ to go to infinity, as long as \reff{cond:H} remains valid (see Section \ref{ssect:Drift}). 
With this assumption, we can prove the cut-off behavior of $\T{a}{\ze}$ (Proposition \ref{prop:1}) and that the typical time scale $\EE{\T{a}{\ze}}$ 
of this convergence is negligible compare to the mean escape time from zero to $a$ (Proposition \ref{prop:2}). If in addition we have a control 
on trajectories outside the branch which contains $a$ (Conditions \reff{hypKK'}), we get the escape behavior of $\T{\ze}{a}$ (Theorem \ref{theo:1}).

Note that until the drift condition is satisfied, our results apply to any tree. 
No particular assumptions are made on the degrees which can be non-homogeneous and non-bounded.\\

The paper is organized as follows. 
We first recall the general definition of the two types of behaviors in Section \ref{ssect:def}. 
Then, we define our birth-and-death chains model on trees (Sections \ref{ssect:Tree}, \ref{ssect:Chain} and \ref{ssect:HT}), 
and the strong drift toward the root in Section \ref{ssect:Drift}. The main results are given in Section \ref{ssect:results} and we discuss 
the basic example of regular trees in Section \ref{ssect:example}. 
Section \ref{sect:mht} contains exact expressions for the mean hitting times and their second moments. 
Finally, the main results are proven in Section \ref{sect:proofmain}.

\subsection{Cut-off and escape behaviors}\label{ssect:def}
Cut-off and escape behaviors, will be studied at the level of hitting times.  
Both types of behavior are asymptotic, in the sense that they are characterized by what happens when a certain parameter $a$ diverges.  
Let us recall here the relevant definitions.
\begin{definition}\label{DefCut/Esc}
\begin{itemize}
\item[(i)] A family of random variables $U^{(a)}$ \emph{exhibits cut-off behavior at mean times} if 
\begin{equation}
\frac{U^{(a)}}{\E{U^{(a)}}} \;\tend{a}{\infty}{\rm Proba}\; 1\;.
\label{eq:0.1}
\end{equation}
[equivalently, $\lim_{a\to\infty}\mathbb{P}\bigl(U^{(a)} > c\, \mathbb{E}[U^{(a)}]\bigr)= 1$ for $c<1$ and $0$ for $c>1$].
\item[(ii)] A family of random variables $V^{(a)}$  \emph{exhibits escape-time behavior at mean times} if 
\begin{equation}
\frac{V^{(a)}}{ \E{V^{(a)}}} \;\tend{a}{\infty}{\mathcal{L}}\; \exp(1) \;.
\label{eq:0.2}
\end{equation}
\end{itemize}
\end{definition}

We refer to \cite{barberfer09} for a discussion about the motivations of these definitions, as well as general sufficient conditions for them to occur.


\section{Model and results}\label{sect:Model}
In order to study cut-off and escape behaviors, we will consider families of random walks $\Xa t$ defined on a tree $I_a$. 
Each of the chains we are defining below is the discrete time counterpart of the birth-and-death processes on trees studied in~\cite{YcartMart}.
\subsection{The trees}\label{ssect:Tree}
\begin{figure}[htbp]
\centering
\scalebox{0.5}{\input{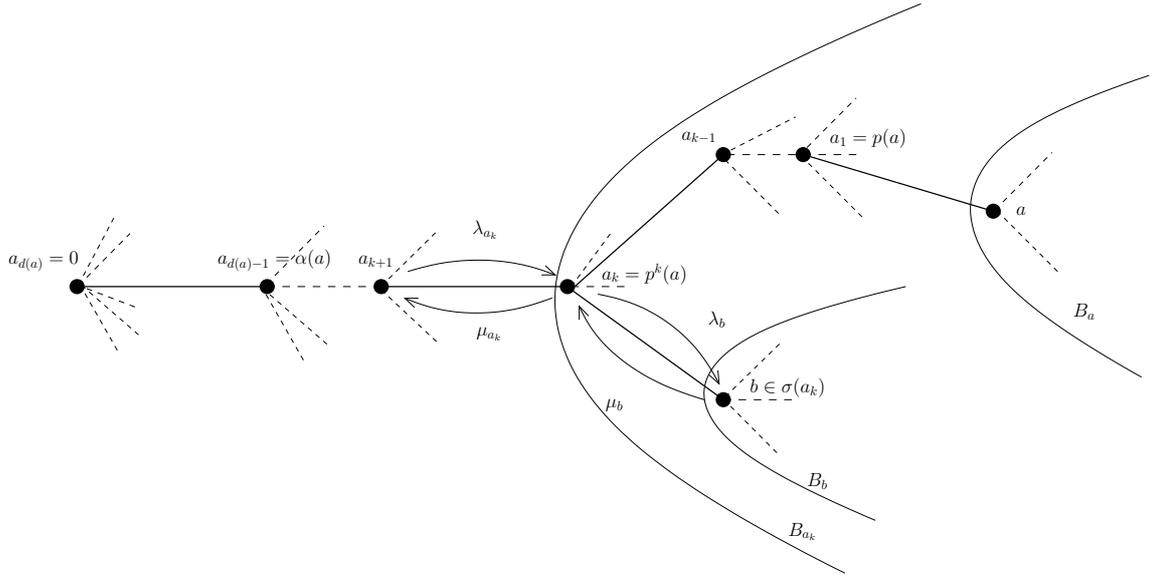}}
\caption{The path from $a$ to zero: sites, transitions and branches}\label{fig:1}
\end{figure}

A tree is an undirected connected graph $G=(I,E)$ with no ``nontrivial closed loop'', where $I$ is the set of vertices (or nodes) 
and $E$ the set of edges (or links).\\
We define a partial order $\preceq$ on $I$ by choosing a node (denoted by $0$) to be the root of the tree:
for $x,y \in I$, we say that $x$ is before $y$  ($x\preceq y$) if $x=y$ or if it exists a path from $y$ to zero containing~$x$.\\
Thus each element $x$ of the tree (except the root) has a unique parent denoted by $\p x$:
\begin{equation}
\forall x \in J\,,\;\exists!\; \p x\; s.t.\quad \p x\preceq x \quad and \quad (\p x,x)\in E\,,
\end{equation}
where $J\bydef I\setminus\{0\}$.\\
There also exists a unique path from $x$ to zero, and we denote by $d(x)$ its length (depth) and by $\ell(x)$ the set of vertices 
on that path (except $0$):
\begin{equation}\label{def:la}
\forall x \in J\;,\;\ell(x)\,\bydef\,\{x_0, x_1,\ldots,x_{d(x)-1}\}\,,
\end{equation}
where $x_i\bydef\gp i x$ for $i=0,\ldots,d(x)-1$. Note that according to this notation $x_0=x$ and~$x_{d(x)}=0$. 
For the sake of notations, we denote by $\alf x$ the last state before zero on that path: $\alf x=x_{d(x)-1}$.\\
Let us also define $\sig x$ the set of the children of $x$, and $\B x$ the branch stemming from $x$: 
\begin{equation}
 \sig x\,\bydef\{y:x=\p y\} \,,\, \text{ and } \; \B x\bydef\{y:x\preceq y\} \,.
\end{equation}
In Section \ref{sect:mht}, we will make use of the following decompositions, for any $k=0,\ldots,d(a)$, of the branch $\Ba k$ stemming from the $k$-th 
parent $a_k=\gp k a$ of some $a\in J$, and its complementary $\Babar k$ (see Figure \ref{fig:1}):
\begin{equation}\label{DecompBa}
 \Ba k \,=\,\B a\,\cup\,\biggl(\,\bigcup_{l=1}^k\,\Bigl( \{a_l\}\,\bigcup\limits_{\substack{c\in\siga l\\  c\neq a_{l-1}}}\,\B c\,\Bigr)\,\biggr)\,,
\end{equation}
\begin{equation}\label{DecompBabar}
 \Babar k \,=\,\bigcup_{l=k+1}^{d(a)}\,\Bigl( \{a_l\}\,\bigcup\limits_{\substack{c\in\siga l\\  c\neq a_{l-1}}}\,\B c\,\Bigr)\,.
\end{equation}
We also need to define the set $\Ca k$
\begin{equation}\label{DecompCa}
 \Ca k \,\bydef\,\bigcup_{l=k+1}^{d(a)}\,\Bigl( \,\bigcup\limits_{\substack{c\in\siga l\\  c\neq a_{l-1}}}\,\B c\,\Bigr)\,.
\end{equation}
Note that $\Ca k$ is the complementary of the path from $a_k$ to zero in $\Babar k$: 
$\displaystyle\Ca k=\Babar k\setminus\bigcup_{l=k+1}^{d(a)}\{a_l\}\,$.


\subsection{The chains}\label{ssect:Chain}
We consider irreducible discrete time birth and death chains $\X t$ on $I$ defined by transition probabilities $\l x$ ``from parent to children'' 
and $\m x$ ``from children to parent'' for each $x\in J$.\\
That is $\l x\bydef P(\p x,x)$ and $\m x\bydef P(x,\p x)$, where $P(x,y)\bydef \prr{\X{t+1}=y\mid\X t=x}$.\\
The $\l i$'s and $\m i$'s sum to $1$ in each site: 
\begin{equation}\label{sommeprobasite}
 \m x+\sum_{y\in\sig x} \l y\,=\,1\,,\;\;\forall x\in J \,,\;\text{ and }\,\sum_{x\in\sig 0} \l x\,=\,1\,,
\end{equation}
and since the chain is irreducible they are supposed to be non-null. \\
Note that, as we will consider families of chains $\Xa t$, the $\lambda$'s and $\mu$'s may vary with $a$.\\
As it doesn't change our results and computations, we can also add (at least one) non-null ``waiting probabilities'' $\kappa_x\bydef P(x,x)$ 
in order to make the chain aperiodic. In that case, equations \reff{sommeprobasite} become $\m x+\kappa_x+\sum_{y\in\sig x} \l y=1$, for all $x\in J$, 
and $\kappa_{\ze}+\sum_{x\in\sig 0} \l x=1$.\\
If the tree is finite, then the chain is positive recurrent and admits a unique reversible invariant probability measure (i.p.m.) $\pi$, 
which can be expressed for each $x\in J$ in terms of the transition probabilities
\begin{equation}\label{invmes}
 \pi(x)\,=\,\pi(0)\,\prod_{y\in\ell(x)}\,\frac{\l y}{\m y}\,.
\end{equation}
This expression is a consequence of the reversibility of $\pi$, which ensures that for each $x\in J$ 
we have $\piii{\p x} \, \l x \,= \,\pii x \, \m x$. Then, iteration along the path $\ell(x)$ from $x$ to zero gives \reff{invmes}, 
where the normalization constant $\pii0$ is chosen such that $\sum_{x\in I}\,\pii x\,=\,1$, that is:
\begin{equation}
 \pii0\,=\,\frac{1}{\sum_{x\in I}\,\prod_{y\in\ell(x)}\,\frac{\l y}{\m y}} \;.
\end{equation}
In the case where the tree is not finite, we have to suppose that the quantity $\sum_{x\in I}\,\prod_{y\in\ell(x)}\,\frac{\l y}{\m y}$ is finite,
 and then the chain is positive recurrent and \reff{invmes} gives the unique i.p.m. $\pi$.

\subsection{Hitting times}\label{ssect:HT}
As mentioned above cut-off and escape behavior  are studied at the level of hitting times.\\ 
More precisely, the families of random variables under consideration are the hitting times (also called first passage times in the literature) 
of zero starting from some $a\in J$ (resp. times from zero to $a$), denoted by $\T{a}{\ze}$  (resp.$\T{\ze}{a}$), with
\begin{equation}
\T{x}{y} \;=\; \min\bigl\{ t\ge 0: \x{x}{t} =y\bigr\}\;,
\end{equation}
where the upper index $x$ means that this is the initial state of the chain ($\x{x}{0}=x$).\\
Our results rely on exact expressions for the first and second moments of $\T{\gp j a}{\gp n a}$, for any $0\leq j,n\leq d(a)$, which will be given 
in Section \ref{sect:mht}.\\
Let us just mention now, for all $x\in J$, the expression in terms of the invariant measure for the mean time $\E{\T{x}{\p x}}$ from child to parent 
(see Lemma \ref{lem:Txpx}):
\begin{equation}\label{ta}
 \E{\T{x}{\p x}}\,=\,\frac{\pii{\B x}}{\m x\,\pii x} \;.
\end{equation}
This expression is a key ingredient in our computations 
and allows also to interpret our drift condition as an exponential decay of $\pi$ (see next section). 

\subsection{The drift condition}\label{ssect:Drift}
As seen in Section \ref{ssect:def}, cut-off and escape behaviors are defined for families of variables indexed by some diverging parameter. 
In the present work, we take a state $a\in J$ such that the length $d(a)$ of the path $\ella$ from $a$ to zero tends to infinity 
to be our asymptotical parameter. \\
We denote by $a\to\infty$ the limit when $d(a)\to\infty$, and we define by $\ell$ the path ``from zero to infinity'' 
followed by $a$ when $d(a)\to\infty$: $\displaystyle\ell\,=\,\lim_{a\to\infty}\,\ella$. \\
Let us give now the definition of strong drift which applies to families of previous defined chains $\Xa t$ on trees $I_a$.
When it is necessary, we mark with an upper index $(a)$ the quantities corresponding to $\Xa t$.

\begin{definition}\label{def:drift}
The family of birth-and-death chains $\Xa t$ has a \emph{strong drift toward 0 along the path $\ell$} (we note $0^{\ell}$-SD) if  
\begin{itemize}
\item[(i)] the transition probabilities $\mua x$ satisfy
\begin{equation} \label{cond:Kmu}
\inf_{b\in\B{\alf a}} \mua b \defby K_{\mu} > 0 \;,\;\text{for all}\;a\in\ell \,,
\end{equation}
\item[(ii)] the constant
\begin{equation} \label{def:Ka}
K_a\;\bydef\; \sup_{b\in\B{\alf a}}\, \mua b\,\EE{\Ta{b}{\p b}} \,,
\end{equation} 
satisfies 
\begin{equation} \label{cond:H}
 \frac{K_a^2}{\EE{\Ta{a}{\ze}}}\linf{a}0 \;.
\end{equation}
\end{itemize}
\end{definition}

Condition $(i)$ ensures that transitions to the left remain non null asymptotically in $a$. It also provides a useful bound (see Lemma \ref{lem:RQK}). 
The drift is given by $(ii)$ and can be interpreted in two ways. Note first that by \reff{ta}, $K_a$ can be written
\begin{equation}\label{eq:Ka}
 K_a \,=\, \sup_{b\in\B{\alf a}}\, \frac{\piia{\B b}}{\piia b} \,.
\end{equation}
Thus, condition $(ii)$ can be seen as an exponential decay of the i.p.m. $\pi^{(a)}$ along the branch $\B{\alf a}$ which contains $a$:
we have $\piia b \,\leq\, e^{-\gamma_a\,d(b)}$, with $\gamma_a\,=\,-\log(1-1/K_a)$, valid for all $b\in\B{\alf a}$. 
That is for any state of any subbranch of the branch containing $a$, since by definition $\B{\alf a}=\bigcup_{x\in\ella}\B x$. 
In terms of energy profile, $\gamma_a$ also provides a lower bound to the energy profile of any subbranch of the branch containing $a$.
We refer to \cite{barberfer09} for a more detailed discussion about these interpretations of our drift condition.

\subsection{Results}\label{ssect:results}

For the fixed path $\ell$ which contains the sequence of $a$'s going to infinity, the $0^{\ell}$-SD drift condition is a sufficient condition for the 
cut-off behavior of hitting times $\Ta{a}{\ze}$:

\begin{prop}\label{prop:1}
If the family $\Xa t$ of irreducible discrete time birth and death chains on $I_a$ satisfies the $0^{\ell}$-SD condition, then
$\displaystyle \lim_{a\to\infty}\,\var \bigl{(}\Ta{a}{\mbox{\tiny\em0}}\,/\,\EE{\Ta{a}{\mbox{\tiny\em0}}}\bigr{)}\,=\,0\,$.\\
As a consequence, the random variables $\Ta{a}{\mbox{\tiny\em0}}$ exhibit cut-off behavior at mean times.
\end{prop}

Moreover, under this drift condition the typical time-scale $\EE{\Ta{a}{\ze}}$ of the abrupt convergence toward zero is negligible comparing 
to the mean times $\EE{\Ta{\ze}{a}}$ of the escape from zero to $a$:

\begin{prop}\label{prop:2}
If the family $\Xa t$ of irreducible discrete time birth and death chains on $I_a$ satisfies the $0^{\ell}$-SD condition, then
 $\displaystyle\lim_{a\to\infty}\EE{\Ta{a}{\mbox{\tiny\em0}}}\,/\,\EE{\Ta{\mbox{\tiny\em0}}{a}}\,=\, 0\,$.
\end{prop}

This time-scales property is the key element of the presence of escape behavior for $\Ta{\ze}{a}$. But, as seen in Section 3 of \cite{barberfer09}, 
this requires also that starting from any state of the state space, the mean hitting times of zero are negligible comparing to $\EE{\Ta{\ze}{a}}$, 
and that the sequence $\Ta{\ze}{a}/\EE{\Ta{\ze}{a}}$ is uniformly integrable.
The requirement that $\sup_{x\in I_a}\EE{\Ta x\ze}/\EE{\Ta\ze a} \linf a 0$ is obtained if we suppose, in addition to the $0^{\ell}$-SD condition, that 
the mean hitting times of zero are comparable. That is if the ratio $\sup_{x\in I_a}\,\EE{\Ta{x}{\ze}}/\EE{\Ta{a}{\ze}}$ remains bounded for all $a$. 
For the uniform integrability of $\Ta{\ze}{a}/\EE{\Ta{\ze}{a}}$, we need to control $\EE{\T{\ze}{a}^2}$ (see Corollary \ref{coro:maj_T2:parent-enfant}). 
A sufficient condition is to have $\displaystyle \inf_{b\in I_a} \, \mua b>0$ uniformly in $a$, and ${K'_a}^2 / \EE{\Ta{\ze}{a}}\linf{a} 0 \,$, where 
\begin{equation} \label{def:K'a}
K'_a\,\bydef\, \sup_{b\in I_a}\, \mua b\,\EE{\Ta{b}{\p b}} \,=\, \sup_{b\in I_a}\, \frac{\piia{\B b}}{\piia b} \;.
\end{equation} 
As above, the second condition is obtained if in addition to the $0^{\ell}$-SD condition, the ratio ${K'_a}^2 / \EE{\Ta{a}{\ze}}$ remains bounded for all $a$
(see Section \ref{ssect:prooftheo:1}).

\begin{theo}\label{theo:1}
 If the family $\Xa t$ of irreducible discrete time birth and death chains on $I_a$ satisfies the $0^{\ell}$-SD condition, 
 and if there exist $K$, $K'$ and $K'_{\mu}$ such that 
\begin{equation}\label{hypKK'}
 \frac{\sup_{x\in I_a}\,\EE{\Ta{x}{\mbox{\tiny\em0}}}}{\EE{\Ta{a}{\mbox{\tiny\em0}}} } \,\leq\,K \,,
\quad  \frac{{K'_a}^2}{\EE{\Ta{a}{\mbox{\tiny\em0}}} } \,\leq\,K' \,,
\,\text{ and }\quad  \inf_{b\in I_a} \mua b \defby K'_{\mu} > 0 \quad \forall a\in \ell\,,
\end{equation}
then:
 \begin{enumerate}
  \item the random variables $\Ta{a}{\mbox{\tiny\em0}}$ exhibit cut-off behavior at mean times.
  \item the random variables $\Ta{\mbox{\tiny\em0}}{a}$ exhibit escape-time behavior at mean times.
 \end{enumerate}
\end{theo}

Moreover, as the quantity $K_a$ is a supremum taken over all sites of the branch $\B{\alf a}$ which contains $a$, these assumptions are sufficient 
to have cut-off and escape behaviors between zero and every state $b$ of $\B{\alf a}$ for which $\EE{\Ta{b}{\ze}}$ diverges with $\EE{\Ta{a}{\ze}}$. 
Let us define
\begin{equation}
 L(a)\,\bydef\,\{\, b\in\B{\alf a}: \, \exists \, C \;\text{s.t.}\; \EE{\Ta{b}{\mbox{\tiny\em0}}}\geq C \, \EE{\Ta{a}{\mbox{\tiny\em0}}} \,\}\,. 
\end{equation}
\begin{coro}\label{coro:1}
Theorem \ref{theo:1} holds for $\Ta{b}{\mbox{\tiny\em0}}$ and $\Ta{\mbox{\tiny\em0}}{b}$ for all $b\in L(a)$.
\end{coro}

As we will see in Section \ref{sect:proofmain} (Lemmas \ref{lem:cut-off} and \ref{lem:escape}), the conditions needed to obtain the two 
time scales $\Bigl(\,\EE{\Ta{a}{\ze}}\,/\,\EE{\Ta{\ze}{a}}\linf{a} 0\,\Bigr)$ -- which are typical of escape behavior -- are actually weaker than those 
which imply cut-off behavior. 
But unlike the case when the chain is defined on $\mathbb N$ where the drift condition implies both cut-off and escape behaviors 
(see \cite{berbarfer08}, Theorem 2.2), here the $0^{\ell}$-SD condition implies only the cut-off behavior directly. 
As the chain started in zero could escape to other branches than $\B{\alf a}$, we need in addition to have a ``control'' on that part of the 
trajectories from zero to $a$ (conditions \reff{hypKK'}) to get the escape behavior of $\Ta{\ze}{a}$.\\
The assumption that $\sup_{x\in I_a}\,\EE{\Ta{x}{\ze}}/\EE{\Ta{\ze}{a}} \linf a 0$ is the analogue of conditions (2.17) in Theorem~1 of 
\cite{barberfer09} for chains on $\mathbb Z$. This is a central argument of the proof of Theorem 2 in \cite{barberfer09} which gives sufficient 
conditions for escape behavior of $\Ta{\ze}{a}$. The idea is that there exists an intermediate time scale $\Delta_a$ between $\EE{\Ta{a}{\ze}}$ and 
$\EE{\Ta{\ze}{a}}$ such that the walk will almost surely (in the limit $a\to \infty$) visit zero in an interval of length $\Delta_a$.
Thus, the walk started in zero will make a great number of unfructuous attempts to escape the attraction of zero before eventually hitting $a$. 
The exponential character (randomness without memory) of $\Ta{\ze}{a}$ is a consequence of this situation: each trial is followed by a return to zero, 
and then by Markovianness, the walk starts again in the same situation.\\
\begin{rem}
If the tree $I_a$ is reduced to the branch $\B{\alf a}$ which contains $a$, 
the conditions of Theorem \ref{theo:1} reduce to $0^{\ell}$-SD plus $\displaystyle\sup_{x\in \B{\alf a}}\,\EE{\Ta{x}{\ze}}/\EE{\Ta{a}{\ze}}$ 
bounded for all $a$.
\end{rem}
\begin{rem}
 Propositions \ref{prop:1} and \ref{prop:2} above also apply 
when the chains  $\Xa t$ are defined on an infinite tree $I$ (i.e. when $I_a=I$ infinite for all $a$). We can also apply Theorem \ref{theo:1} 
in that case by considering the restrictions of the chains to $I_a\bydef \{x\in I: d(x)\leq d(a)\}$.
\end{rem}

Another feature of the complementarity between both phenomena, is that the successful final excursion which leads to $a$ is ``cut-off like'' 
(see Theorem 1 and Section 7 in \cite{barberfer09}): if we consider the hitting times $\Tat{\ze}{a}$ of $a$ after the last visit to zero, 
and  $\Tat{a}{\ze}$ (hitting time of zero after last visit to $a$), we can show by a reversibility argument 
that they have the same laws. Thus, cut-off and escape behaviors are related by time inversion.

\subsection{Example: biased random walk on a regular tree}\label{ssect:example}
A regular tree is a tree $I$ in which each node has the same number $r$ of children. In its infinite version (also called Bethe lattice), 
the degree (number of neighbors) of each node is  $r+1$, except the root which is of degree $r$. When we consider the finite version 
$I_a\bydef \{x\in I: d(x)\leq d(a)\}$ for some $a$, the boundary nodes (sites $x$ s.t. $d(x)=d(a)$) -- also called leaves of the tree -- 
are dead ends (of degree 1). In that case, we talk about a Cayley tree. We refer to \cite{DorGovMend08,DoroLectures10} for discussions 
about the differences between Bethe lattice and Cayley tree, and to Chapter 5 of \cite{AldousFill} for the biased random walk on a regular tree.\\

We consider biased random walks $\Xa t$ on $I_a$ which transitions are defined by
\begin{equation}
 \l x\bydef\,\frac{1}{\lambda+r}\;\;\text{and }\;\; \m x\bydef\,\frac{\lambda}{\lambda+r}\;,
\end{equation}
with $0<\lambda<\infty$, for each site $x\neq0$. As transitions must sum to one in each site, we add also ``waiting probabilities'' at root and leaves: 
$\kappa_0\bydef\lambda/(\lambda+r)$ and $\kappa_x\bydef r/(\lambda+r)$ when $d(x)=d(a)$. 
The parameter $\lambda$ represents the bias of the random walk: at each step, the chain moves toward the root of the tree  
with probability $\lambda/(\lambda+r)$ and in the opposite direction with probability $r/(\lambda+r)$. Thus, there is a concentration phenomenon 
near the root level when $\lambda>r$ (localization), and near the leaves level when $\lambda<r$ (delocalization). \\

Due to the regular structure of the tree, we can compute the quantities involved in the $0^{\ell}$-SD condition, and determine when it is satisfied 
or not. Let us start with the invariant measure $\pi$. By \reff{invmes}, for all $x$ we have
\begin{equation}
 \pi(x)\,=\,\pi(0)\,\prod_{y\in\ell(x)}\,\frac{\l y}{\m y}
	\,=\,\pi(0)\,\prod_{y\in\ell(x)}\,\frac{1}{\lambda+r}\frac{\lambda+r}{\lambda}\,=\,\pi(0)\,\frac{1}{\lambda^{d(x)}} \;,
\end{equation}
with $\displaystyle\pi(0)\,=\,1\,/\,\sum_{x\in I_a} \frac{1}{\lambda^{d(x)}} \,$.\\

\begin{rem}\label{rem:MRT}
The computations of mean return times using Kac formula 
$\bigl(\, \EE{\TR{x}}=1/\pi(x)$, where $\TR{x}\bydef\inf\{t>0:\x{x}{t} =x\}$ is the first return time to $x \,\bigr)$
illustrates the localization phenomenon for the biased random walk on a regular tree:
\begin{itemize}
\item localization is related to positive recurrence of the walk when it takes place on an infinite tree. Lyons showed that the biased walk is 
recurrent (return times to any state almost surely finite) when the bias $\lambda$ exceeds the branching number of the tree (here  $\lambda>r$) and 
transient when $\lambda<r$ (see \cite{Lyons90}).\\ Moreover, when $\lambda>r$, the chain is positive recurrent (finite mean return times to any state). 
Using the Kac formula and since for each $k$ there are $r^k$ sites $y$ in $I$ such that $d(y)=k$, we get
\begin{equation}
 \EE{\TR{x}}\,=\,\lambda^{d(x)}\,\sum_{y\in I} \frac{1}{\lambda^{d(y)}}\,=\,\lambda^{d(x)}\,\sum_{k=0}^{\infty} \left(\frac{r}{\lambda}\right)^{k}
    \,<\, \infty\,,\; \text{ if } \lambda>r\,.
\end{equation}
\item in the finite case, although the walk on $I_a$ is always positive recurrent, the localization phenomenon is characterized by very different 
asymptotical behaviors of mean return times to zero and $a$ respectively. When $\lambda>r$ the mean return time to zero remains bounded as $a$ goes 
to infinity (by Kac formula, we get $\EE{\TRa{\ze}}\sim\lambda/(\lambda-r)$ for large $a$). In opposition, the mean return time to $a$ diverges 
with $a$ ($\EE{\TRa{a}}\sim\lambda^{d(a)+1}/(\lambda-r)$ for large $a$).
\end{itemize}
\end{rem} 
\noindent Let us now compute the key quantity $\displaystyle\frac{\pii{\B x}}{\pii x}$ which appears in both expressions of $K_a$ and $\EE{\Ta{a}{\ze}}$ 
(see \reff{eq:Ka}, and Proposition \ref{prop:esp} below). We write
\begin{equation}
 \frac{\pii{\B x}}{\pii x}\,=\,\sum_{b\in\B x}\,\frac{\pii{b}}{\pii x} \,=\,\sum_{b\in\B x}\,\frac{\lambda^{d(x)}}{\lambda^{d(b)}} 
  \,=\,\sum_{k=0}^{d(a)-d(x)}\,\left(\frac{r}{\lambda}\right)^k\,,
\end{equation}
since there are $r^k$ sites $b$ in $\B x$ such that $d(b)=d(x)+k$.\\

When $\lambda\neq r$, this is equal to $\displaystyle \frac{\lambda}{\lambda-r}\left(1-\left(\frac{r}{\lambda}\right)^{d(a)-d(x)+1}\right)$, and thus
\begin{equation}
  K_a \,\bydef\, \sup_{x\in\ella} \,\sup_{b\in\B x}\, \frac{\pii{\B b}}{\pii b}  
		\,=\,  \frac{\lambda}{\lambda-r}\left(1-\left(\frac{r}{\lambda}\right)^{d(a)}\right) \,,
\end{equation}
and 
\begin{eqnarray}
 \EE{\Ta{a}{\ze}} &=& \sum_{x\in\ella}\,\dfrac{\pii{\B x}}{\m x\,\pii{x}} 
	\,=\,\frac{\lambda+r}{\lambda-r}\,\sum_{x\in\ella} \,\left(1-\left(\frac{r}{\lambda}\right)^{d(a)-d(x)+1}\right) \nonumber\\
  &=& \frac{\lambda+r}{\lambda-r}\,
	    \left(d(a)\,-\,\left(\frac{r}{\lambda}\right)^{d(a)}\,\sum_{k=0}^{d(a)-1} \,\left(\frac{\lambda}{r}\right)^{k}\right) \\
&=& \frac{\lambda+r}{\lambda-r}\,
	    \left( d(a)\,+\,\frac{r}{\lambda-r}\left(\left(\frac{r}{\lambda}\right)^{d(a)}-1\right) \right) \nonumber \,.
\end{eqnarray}
In the localization case ($\lambda>r$), we have $\displaystyle\EE{\Ta{a}{\ze}}\sim\frac{\lambda+r}{\lambda-r}\;d(a)$, 
and $\displaystyle K_a\to K\bydef \frac{\lambda}{\lambda-r} $ as $a\to\infty$.\\ The $0^{\ell}$-SD condition is then satisfied, and 
we have cut-off and escape behavior between the root and any leave of $I_a$, since for symmetry reasons, the above computations are the same for any 
state of $I_a$. This result leads to a better understanding of what happens in that case: since $\Ta{a}{\ze}$ exhibits cut-off behavior at mean times 
$\EE{\Ta{a}{\ze}}\sim d(a)$, the walk started in $a$ will go in an almost deterministic way to zero before eventually returning to $a$. Moreover, 
the typical time scale of cut-off trajectories is negligible comparing to that of escape trajectories (computations give $\EE{\Ta{\ze}{a}}\sim\lambda^{d(a)}$ 
for large $a$, which is of the same order than $\EE{\TRa{a}}$).\\
In the delocalization case ($\lambda<r$), $K_a$ and $\EE{\Ta{a}{\ze}}$ are both exponential of order $(r/\lambda)^{d(a)}$ 
for large $a$, and then there is no strong drift. 
The $0^{\ell}$-SD condition is also not satisfied when $\lambda=r$. Computations in that case show that $K_a=d(a)$ and 
$\EE{\Ta{a}{\ze}}=\bigl((\lambda+r)/2\lambda\bigr)\,d(a)\bigl(d(a)+1\bigr)$ and thus the ratio $K_a^2/\EE{\Ta{a}{\ze}}$ does not converge to zero 
when $a$ goes to $\infty$.\\



\section{Mean hitting times}\label{sect:mht}
We already mentioned (without proving it) in Section \ref{ssect:HT}, that the expression of $\E{\T{x}{\p x}}$ in terms of the invariant measure $\pi$ 
(equation \reff{ta}) is a key ingredient in our computations. 
However, we will also need exact expressions for the first and second moments of hitting times between sites 
along the path from zero to some fixed state $a\in J$. 
We give a list of these expressions in Section \ref{ssect:exact}. We start with first and second moments for $\T{x}{\p x}$ (Lemma \ref{lem:Txpx}), 
and $\T{\p x}{x}$ (Lemma \ref{lem:Tpxx}). Proposition \ref{prop:esp} gives the mean hitting times between sites $a_i\in \ella$ for some $a\in J$, 
and Proposition \ref{prop:moment2} their second moments.
Finally we give an idea of the proofs in Section \ref{ssect:proofexact}.

\subsection{Exact expressions}\label{ssect:exact}

\begin{lem}\label{lem:Txpx}For any $x\in J$
\begin{equation}\label{eq:lem:tb}
 \E{\T{x}{\p x}}\,=\,\frac{\pii{\B x}}{\m x\,\pii x} \;,
\end{equation}
 \begin{equation}\label{eq:lem:ub}
  \EE{\T{x}{\p x}^2}\,=\,\frac{2}{\m x\,\pii{x}} \,\sum_{b\in\B x}\, \frac{\pii{\B b}^2}{\m b\,\pii{b}} \,-\,\EE{\T{x}{\p x}} \;.
 \end{equation}
\end{lem}

\begin{lem}\label{lem:Tpxx}For any $x\in J$
\begin{equation}\label{eq:lem:tbarb}
\E{\T{\p x}{x}}\,=\,\frac{1-\pii{\B x}}{\m x\,\pii x} \;,
\end{equation}
 \begin{equation}\label{eq:lem:ubarb}
  \EE{\T{\p x}{x}^2}\,=\, \frac{2}{\m x\,\pii{x}} 
	    \,\left(\, \sum_{c\in C_x}\, \frac{\pii{\B c}^2}{\m c\,\pii{c}} \,+\,\sum_{b\in\ell(x)}\,\frac{\pii{\Bbar b}^2}{\m b\,\pii{b}} \,\right)
      \,-\,\EE{\T{\p x}{x}} \;.
 \end{equation}
\end{lem}

\begin{prop}\label{prop:esp}
 For $a\in J$, recall that we denote by $a_i$ the $i$-th parent of $a$ ($\, a_i=\gp i a \,$).\\ 
Let $0\leq j<n\leq d(a)$, we have
    \begin{equation}\label{T:enfant-parent}
    \EE{\T{a_j}{a_n}} \,=\, \sum_{k=j}^{n-1}\,\dfrac{\pii{\Ba k}}{\ma k\,\pii{a_k}} \;,
    \end{equation}
     \begin{equation}\label{T:parent-enfant}
    \EE{\T{a_n}{a_j}} \,=\, \sum_{k=j}^{n-1}\,\frac{1-\pii{\Ba k}}{\ma k\,\pii{a_k}} \;.
    \end{equation}
And thus
    \begin{equation}\label{T:parent+enfant}
      \EE{\T{a_j}{a_n}}\,+\,\EE{\T{a_n}{a_j}} \,=\, \sum_{k=j}^{n-1}\,\dfrac{1}{\ma k\,\pii{a_k}} \;.
    \end{equation}
\end{prop}

\begin{prop}\label{prop:moment2}
 For $0\leq j<n\leq d(a)$, we have
\begin{equation}\label{T2:enfant-parent}
    \EE{\T{a_j}{a_n}^2} \,=\, 
	\sum_{k=j}^{n-1}\,\frac{2}{\ma k\,\pii{a_k}}\,\Biggl(\, \sum_{b\in\Ba k}\, \frac{\pii{\B b}^2}{\m b\,\pii{b}}
			\, + \,\pii{\Ba k}\,\EE{\T{a_{k+1}}{a_n}} \,\Biggr) \,-\, \,\EE{\T{a_j}{a_n}} \;,
\end{equation}
and
\begin{equation}\label{T2:parent-enfant}
    \EE{\T{a_n}{a_j}^2} =  \sum_{k=j}^{n-1}\frac{2}{\ma k\pii{a_k}}\Biggl(\,\sum_{c\in\Ca{k}} \frac{\pii{\B c}^2}{\m c\,\pii{c}}
			    \,+ \sum_{b\in\ell(a_{k})} \frac{\pii{\,\Bbar b}^2}{\m b\,\pii{b}}  
			    \,+\, \pii{\Babar k} \EE{\T{a_{k}}{a_j}} \Biggr) - \EE{\T{a_n}{a_j}} \,.
\end{equation}
\end{prop}

\begin{coro}\label{coro:maj_T2:parent-enfant}
For $0\leq j<n\leq d(a)$, we have
\begin{equation}\label{eq:maj_T2:parent-enfant}
  \EE{\T{a_n}{a_j}^2} \,\leq \,\EE{\T{a_n}{a_j}} \,\left(\,2 \,\left(\, \frac{{K'_a}^2}{K'_\mu}\,+\,\EE{\T{\mbox{\tiny\em0}}{a_j}}\,\right)\,-1\,\right)\;.
\end{equation}
\end{coro}
\noindent{\bf{Proof}}\\
This a consequence of equation \reff{T2:parent-enfant}. We first use the fact that $\Bbar b\subseteq \Babar k$ for all $b\in \ell(a_k)$ to write
\begin{equation}
 \sum_{b\in\ell(a_{k})} \frac{\pii{\,\Bbar b}^2}{\m b\,\pii{b}} \,\leq\,\pii{\,\Babar k}\,\sum_{b\in\ell(a_{k})} \frac{\pii{\,\Bbar b}}{\m b\,\pii{b}}
    \,=\,\pii{\,\Babar k}\,\EE{\T{\ze}{a_k}} \,.
\end{equation}
Then, as $\Ca k\subset \Babar k$, and since by definition $\displaystyle K'_a \,=\, \sup_{b\in I_a}\, \frac{\pii{\B b}}{\pii b}$ 
and $\displaystyle K'_\mu \bydef \inf_{b\in I_a} \mua b$, we have
\begin{equation}
 \sum_{c\in\Ca{k}} \frac{\pii{\B c}^2}{\m c\,\pii{c}} \,\leq\, \frac{{K'_a}^2}{K'_\mu}\, \sum_{c\in\Ca{k}} \pii{c} 
\,=\, \frac{{K'_a}^2}{K'_\mu} \, \pii{\Ca{k}} \,\leq\, \frac{{K'_a}^2}{K'_\mu} \, \pii{\Babar k} \,.
\end{equation}
Finally, from equation \reff{T2:parent-enfant}
\begin{equation}
 \EE{\T{a_n}{a_j}^2} \,\leq \, 2\,\sum_{k=j}^{n-1}\,\frac{\pii{\Babar k}}{\ma k\pii{a_k}}
			    \Biggl(\,\frac{{K'_a}^2}{K'_\mu} + \EE{\T{\ze}{a_k}} + \EE{\T{a_{k}}{a_j}}\,\Biggr) - \EE{\T{a_n}{a_j}} \,.
\end{equation}
Hence the result, since $\T{\ze}{a_k}+\T{a_{k}}{a_j}=\T{\ze}{a_j}$ and then using equation \reff{T:parent-enfant}.

\begin{flushright}
 $\Box$ 
\end{flushright}

\subsection{Proofs}\label{ssect:proofexact}
The above formulas can be proven via the standard method of difference equations which rely on the following decomposition of $\EE{F(\T{a_k}{a_n})}$ 
for any function $F$:

\begin{lem}\label{lem:claim}
For all $0<k<d(a)$, $k\neq n$, we have
\begin{equation}\label{eq:claim}
  \EE{F(\T{a_k}{a_n})} \,=\, \ma k\, \EE{F(\T{a_{k+1}}{a_n}+1)}\,+\,\sum_{c\in\siga k}\l c \,\EE{F(\T{c}{a_n}+1)} \;.
 \end{equation}
\end{lem}

\noindent{\bf{Proof}}\\
This is obtained by decomposing according to the first step of the chain
\begin{equation}
 \EE{F(\T{a_k}{a_n})} = \ma k\, \EE{F(\T{a_k}{a_n})\Bigm|X(1)=a_{k+1}}\,+\sum_{c\in\siga k} \l c \,\EE{F(\T{a_k}{a_n})\Bigm|X(1)=c}\,,
\end{equation}
and by Markovianness, for $x=a_{k+1}$ and $c\in\siga k$, we have
\begin{eqnarray}
\EEE{F(\T{a_k}{a_n})\Bigm|X(1)=x} &=& \sum_{l\ge 1} F(l)\, \prr{\T{a_k}{a_n}=l\Bigm|X(1)=x}\nonumber\\
&=& \sum_{l\ge 1} F(l)\, \pr{\T{x}{a_n}=l-1}\nonumber\\
&=& \sum_{l\ge 0} F(l+1)\, \pr{\T{x}{a_n}=l}\nonumber\\
&=& \EEE{F(\T{x}{a_n}+1)}\;.
\end{eqnarray}
Hence the result.
\begin{flushright}
 $\Box$ 
\end{flushright}

Together with Lemma \ref{lem:Txpx} which is the discrete analogue of Lemma 1 of \cite{YcartMart} for continuous chains, Lemma \ref{lem:claim} 
allows us to prove Propositions \ref{prop:esp} and \ref{prop:moment2}: we fix $0\leq n\leq d(a)$ and apply equation \reff{eq:claim} with $F(x)=x$ 
(resp. $F(x)=x^2$) to obtain recursive equations for $\EE{\T{a_k}{a_n}}$, (resp. $\EE{\T{a_k}{a_n}^2}$) with $0<k<d(a)$. Then we can iterate these 
equations in both directions up to boundary conditions -- which are particular cases of \reff{eq:claim} for $k=0$ and $d(a)$ -- 
and solve them using reversibility and equations \reff{eq:lem:tb} and \reff{eq:lem:ub}.\\
Another way of proving both propositions is to use, together with Lemmas \ref{lem:Txpx} and \ref{lem:Tpxx}, 
the fact that $\T{a_j}{a_n}$ (resp. $\T{a_n}{a_j}$) is the sum of independant $\T{x}{\p x}$ (resp. $\T{\p x}{x}$):
\begin{equation}\label{eq:SumIndepT}
 \T{a_j}{a_n}\,=\,\sum_{k=j}^{n-1}\,\T{a_k}{\pa k} \quad\text{and } \quad \T{a_n}{a_j}\,=\,\sum_{k=j}^{n-1}\,\T{\pa k}{a_k} \;.
\end{equation}

We give now a sketch of those proofs.\\
\noindent{\it{First moments:}}\\
Althought they can be proven directly by difference equations, formulas \reff{eq:lem:tb} and \reff{eq:lem:tbarb} are proven in \cite{Woess} 
(Proposition 9.8) by an alternative method. 
And thus \reff{T:enfant-parent} and \reff{T:parent-enfant} of Proposition \ref{prop:esp} follow from \reff{eq:SumIndepT}.\\

\noindent{\it{Second moments:}}\\ 
The proof of \reff{eq:lem:ub} follows the same lines as the ones of Lemma 1 in \cite{YcartMart}, where the application of equation \reff{eq:claim} gives the 
recursive equation $\m x \,t_x = 1+\sum_{c\in\sig x}\l c \,t_c$ for $t_x=\EE{\T{x}{\p x}}$. 
We now apply \reff{eq:claim} with $F(x)=x^2$, $a_k=x$ and $a_n=\p x$ to obtain the recursive equation for $\EE{\T{x}{\p x}^2}$:
\begin{equation}
 \m x \,\EE{\T{x}{\p x}^2} \,=\, 2\,\m x (\EE{\T{x}{\p x}})^2 \,+\, \sum_{c\in\sig x}\l c \,\EE{\T{c}{\p c}^2} \,-\,1 \;.
\end{equation}
Then we prove by induction that formula \reff{eq:lem:ub} satisfies the above equation. As in \cite{YcartMart}, we have first to suppose that 
the tree is finite, since the first step of induction is to consider a leave $x$ (that is when $\sig x=\varnothing$). For the infinite case, 
we consider the limit $n\to\infty$ of restrictions to the truncated trees $I_n\bydef \{x\in I: d(x)\leq n\}$ of level $n$.\\
For equation \reff{eq:lem:ubarb}, we prove by induction that it satisfies the recursive equation obtained from \reff{eq:claim} 
with $F(x)=x^2$, $a_k=\p x$ and $a_n=x$. We don't need to restrict ourselves to finite trees in that case, since the first step of induction is to 
consider $x$ such that $\p x=0$.\\
Finally, \reff{T2:enfant-parent} and \reff{T2:parent-enfant} of Proposition \ref{prop:moment2} follow from \reff{eq:SumIndepT}. 




\section{Proofs of the main results}\label{sect:proofmain}
The quantities involved in this section are defined with respect to the chain $\Xa t$ on $I_a$, but for the sake of notations we will omit 
the upper index $(a)$. We will prove Propositions \ref{prop:1} and \ref{prop:2} through Lemmas \ref{lem:cut-off} and \ref{lem:escape} 
which are a little bit stronger. They involve the quantities
\begin{equation}\label{def:Qa}
\Q x \bydef \frac{1}{\pii{\B x}} \sum_{b\in\B x}\, \frac{\pii{\B b}^2}{\m b\,\pii b}
\quad;\quad \qa \bydef \sup_{x\in\ella} \Q x\;,
\end{equation}
and
\begin{equation}\label{def:Ra}
 \ra \bydef  \sum_{b\in\ella}\, \frac{\pii{\B b}^2}{\m b\,\pii b} \,,
\end{equation}
which for all $a$ satisfy the inequalities
\begin{lem}\label{lem:RQK}
\begin{equation}\label{ineg:RQK}
 \ra \,\leq\,  \qa \,\leq\,K_a^2\,/\,K_\mu  \,.
\end{equation}
\end{lem}
\noindent{\bf{Proof of Lemma \ref{lem:RQK}}}\\
As $\ella\subset\Ba{d(a)-1}$, we have
\begin{equation}
 \ra \,\leq\,  \sum_{b\in\Ba{d(a)-1}}\, \frac{\pii{\B b}^2}{\m b\,\pii b} \,=\, \pii{\Ba{d(a)-1}}\,\Q{a_{d(a)-1}}  \,,
\end{equation}
and thus $ \ra \,\leq\,  \qa$, since $\pii{\Ba{d(a)-1}}\,\Q{a_{d(a)-1}} \,\leq\,\Q{a_{d(a)-1}} \,\leq\,\qa$.\\
For the right-hand side of \reff{ineg:RQK}, we write
\begin{equation}
 \frac{1}{\pii{\B x}} \, \sum_{b\in\B x} \, \frac{\pii{\B b}^2}{\m b\,\pii b} 
      \,\leq\, \frac{1}{\pii{\B x}} \, \frac{K(x)^2}{K_\mu} \, \, \sum_{b\in\B x}\, \pii{b}
      \,=\,\frac{K(x)^2}{K_\mu} \,,
\end{equation}
where $ K(x)=\sup_{b\in\B x}\left(\pii{\B b}\,/\,\pii b\right) $. \\
We conclude by taking the supremum over $x\in \ella$, since by \reff{eq:Ka}, $\displaystyle K_a=\sup_{x\in\ella}K(x)$.
\begin{flushright}
 $\Box$ 
\end{flushright}
\subsection{Proof of Proposition \ref{prop:1}}\label{ssect:proofprop1}
We recall that 
\begin{equation}
\lim_{a\to\infty} \var\Bigl{(}\frac{\T a\ze}{\EE{\T a\ze}} \Bigr{)} =0
\end{equation}
is a sufficient condition to prove the cut-off behavior of the variables $\T a\ze$ (see \cite{barberfer09}, Section 3).\\
Since $\qa\leq K_a^2\,/\,K_\mu $ (see Lemma \ref{lem:RQK} above), and since by the $0^{\ell}$-SD hypothesis the ratio $K_a^2/\EE{\T a\ze}$ tends to 
zero, the following result implies Proposition \ref{prop:1}.
\begin{lem}\label{lem:cut-off}
 If \begin{equation} \label{cond:Hbis} 
 \frac{\qa}{\EE{\T{a}{\mbox{\tiny\em0}}}}\linf{a}0 \;,
\end{equation}
then $\displaystyle \var \bigl{(}\T{a}{\mbox{\tiny\em0}}/\EE{\T{a}{\mbox{\tiny\em0}}}\bigr{)}\linf a0 \,$. \\
\end{lem}

\noindent{\bf{Proof of Lemma \ref{lem:cut-off}}}\\
For all $a\in J$, we have
\begin{eqnarray}
 \var\Bigl{(}\T a\ze \Bigr{)}&=& \sum_{x\in\ella}\, \var\Bigl{(}\T{x}{\p x} \Bigr{)} \nonumber\\
	&\leq&  \sum_{x\in\ella}\,\EE{\T{x}{\p x}^2}\\
	&\leq&  \sum_{x\in\ella}\,\frac{2}{\m x\,\pii{x}} \,\sum_{b\in\B x}\, \frac{\pii{\B b}^2}{\m b\,\pii{b}} \;, \nonumber
\end{eqnarray}
where the second inequality comes from equation \reff{eq:lem:ub} of Lemma \ref{lem:Txpx}.\\
Thus
\begin{equation}
 \var\Bigl{(}\T a\ze \Bigr{)}\,\leq\, \sum_{x\in\ella}\,\frac{2\,\pii{\B x}}{\m x\,\pii{x}} \; \Q x
      \,\leq \,  2\,\qa\,\sum_{x\in\ella}\,\frac{\pii{\B x}}{\m x\,\pii{x}}\;,
\end{equation}
and by \reff{eq:lem:tb}
\begin{equation}
 \var\Bigl{(}\frac{\T a\ze}{\EE{\T a\ze}}\Bigr{)} \,\leq\, \frac{2\,\qa}{\EE{\T{a}{\ze}}}\,\linf a0 \;.
\end{equation}
\begin{flushright}
 $\Box$ 
\end{flushright}




\subsection{Proof of Proposition \ref{prop:2}}\label{ssect:proofprop2}
As above, by Lemma \ref{lem:RQK} ($ \ra \,\leq\, K_a^2\,/\,K_\mu$) together with the $0^{\ell}$-SD condition, the following lemma implies 
Proposition \ref{prop:2}.
\begin{lem}\label{lem:escape}
 If \begin{equation} \label{cond:Hter} 
 \frac{\ra}{\EE{\T{a}{\mbox{\tiny\em0}}}}\linf{a}0 \;,
\end{equation}
then $\EE{\T{a}{\mbox{\tiny\em0}}}\,/\,\EE{\T{\mbox{\tiny\em0}}{a}}\,\linf{a}0 \,$.\\
\end{lem}
\noindent{\bf{Proof of Lemma \ref{lem:escape}}}\\
It is equivalent to show that
\begin{equation}
 \Ga\,\bydef\,\frac{\EE{\T a\ze}}{\EE{\T a\ze}+\EE{\T\ze a}}\,\linf a0 \;.
\end{equation}
By \reff{T:enfant-parent} and \reff{T:parent+enfant}, we have
\begin{equation}
 \Ga\,=\,\frac{\displaystyle \sum_{x\in\ella}\,\dfrac{\pii{\B x}}{\m x\,\pii{x}} }{\displaystyle \sum_{x\in\ella}\,\dfrac{1}{\m x\,\pii{x}} } \;\;.
\end{equation}
We now consider the probability measure on $\ella$ defined by expectations
\begin{equation} 
\mathcal{E}(F) \;\bydef\; 
    \frac{\displaystyle\sum_{x\in\ella}\,\dfrac{F(x)}{\m x\,\pii{x}} }{ \displaystyle\sum_{x\in\ella}\,\dfrac{1}{\m x\,\pii{x}} }\;\;,
\end{equation}
and the inequality $\mathcal{E}\left(F\right)^2\le \mathcal{E}\left(F^2\right) $ for $F(x)=\pii{\B x}$ implies that 
$\Ga \,\leq\, \ra\,/\,\EE{\T a\ze}$, which concludes the proof.
\begin{flushright}
 $\Box$ 
\end{flushright}


\subsection{Proof of Theorem \ref{theo:1}}\label{ssect:prooftheo:1}

The cut-off behavior of $\T{a}{\ze}$ was already proven in Proposition \ref{prop:1}. 
For the escape behavior of $\T{\ze}{a}$, we need to prove that $\sup_{x\in I_a}\EE{\T x\ze}/\EE{\T \ze a}$ tends to zero and that the sequence 
$\T \ze a\,/\,\EE{\T \ze a}$ is uniformly integrable (see \cite{barberfer09}, Theorem 2).\\ For the first point, we have
\begin{equation}
 \frac{ \sup_{x\in I_a}\,\EE{\T x\ze} }{ \EE{\T\ze a} } 
	\,=\, \frac{ \sup_{x\in I_a}\,\EE{\T x\ze} }{ \EE{\T{a}{\ze}} } \;\; \frac{ \EE{\T{a}{\ze}} }{ \EE{\T\ze a} }
      \,\leq\,K\,\frac{ \EE{\T{a}{\ze}} }{ \EE{\T\ze a} }\linf a 0 \;,
\end{equation}
as a consequence of \reff{hypKK'} and Proposition \ref{prop:2}.\\
For the uniform integrability, we apply Corollary \ref{coro:maj_T2:parent-enfant} with $n=d(a)$ and $j=0$ to get
\begin{equation}
 \frac{\EE{\T{\ze}{a}^2}}{\EE{\T\ze a}^2} \,\leq\, 2+ \frac{2}{K'_{\mu}} \, \frac{{K'_a}^2}{\EE{\T\ze a}} \,\leq\, 2+o(1) \,.
\end{equation}
The second inequality is due to \reff{hypKK'} and Proposition \ref{prop:2}, and thus the sequence $\T\ze a\,/\,\EE{\T\ze a}$ is uniformly integrable.

\subsection{Proof of Corollary \ref{coro:1}}\label{ssect:proofcoro:1}
The key element of this result, is that $K_b$ is equal to $K_a$, as $\B{\alf b}=\B{\alf a}$ for $b\in L(a)$. Then if the $0^{\ell}$-SD condition holds 
for $a$, it also holds for every state $b$ of $L(a)$:
\begin{equation}
 \frac{K_b^2}{\EE{\T{b}{\ze}}} \,=\, \frac{K_a^2}{\EE{\T{b}{\ze}}} \,\leq\, \frac{K_a^2}{C\,\EE{\T{a}{\ze}}} \,\linf a 0 \,.
\end{equation}

Since equation \reff{ineg:RQK} is valid for any state, we can apply Lemma \ref{lem:cut-off} and Proposition \ref{prop:1} holds for all $b$: 
that is $\T{b}{\ze}$ exhibits cut-off behavior for every $b\in L(a)$.\\
It remains to prove the escape behavior of $\T{\ze}{b}$. Let $b\in L(a)$, using \reff{hypKK'} we have
\begin{equation}
 \frac{\sup_{x\in I_a}\EE{\T x\ze}}{\EE{\T\ze b}} 
	\,=\,\frac{\sup_{x\in I_a}\EE{\T x\ze}}{\EE{\T a\ze}} \,\frac{\EE{\T a\ze}}{\EE{\T b \ze}}\, \frac{\EE{\T b\ze}}{\EE{\T \ze b}}
	\,\leq\, \frac{K}{C}\,\frac{\EE{\T b\ze}}{\EE{\T \ze b}} \,,
\end{equation}
and this tends to zero since by Lemma \ref{lem:escape},  Proposition \ref{prop:2} holds for $b$.\\

\noindent Now by Corollary \ref{coro:maj_T2:parent-enfant}, we have
\begin{equation}
 \frac{\EE{\T{\ze}{b}^2}}{\EE{\T\ze b}^2} \,\leq\, 2+ \frac{2}{K_{\mu}} \, \frac{{K'_a}^2}{\EE{\T\ze b}} \,\leq\, 2+o(1) \,,
\end{equation}
since by \reff{hypKK'} and Proposition \ref{prop:2}
\begin{equation}
 \frac{{K'_a}^2}{\EE{\T\ze b}} \,=\,\frac{{K'_a}^2}{\EE{\T a\ze}} \,\frac{\EE{\T a\ze}}{\EE{\T b \ze}}\, \frac{\EE{\T b\ze}}{\EE{\T \ze b}}
      \,\leq\,\frac{K'}{C}\,\frac{\EE{\T b\ze}}{\EE{\T \ze b}} \,\linf a 0 \,.
\end{equation}
Thus $\T \ze b\,/\,\EE{\T \ze b}$ is uniformly square integrable, and $\T{\ze}{b}$ exhibits escape behavior.\\

\noindent {\bf Acknowledgments.} The author wishes to thank an anonymous referee for the careful reading of the manuscript and his helpful 
comments and suggestions to improve the paper.









\newpage
\section*{Appendix: Difference equations for mean hitting times}\label{append}
We give here the detailed proofs of Propositions \ref{prop:esp} and \ref{prop:moment2} of Section \ref{sect:mht} using the standart method of 
difference equations.\\


\noindent{\bf{Proof of Proposition \ref{prop:esp}}}\\
Let $a\in J$ and $0\leq n\leq d(a)$ be fixed, and define $D_x\bydef\,\EE{F(\T{x}{a_n})}$ for all $x\in\ella$. 
Since by \reff{sommeprobasite} the transitions out of each site sum to $1$, equation \reff{eq:claim}, with $F(T)=T$, can be rewritten 
\begin{equation}\label{eq:recc}
\Da k\,=\,\ma k\,\Da{k+1}\,+\, \la{k-1}\,\Da {k-1}\,+\,\sum\limits_{\substack{c\in\siga k\\  c\neq a_{k-1}}}\,\l c \,\D c\,+\,1\,,
\end{equation}
valid for all $0<k<d(a)$, $k\neq n$, with the boundary conditions 
\begin{eqnarray}
\Da 0 \bydef& \EE{\T{a}{a_n}} =& \m a\,\Da1\,+\,\sum_{c\in\sig a}\,\l c \,\D c\,+\,1\,,\\
\Da{d(a)} \bydef& \EE{\T{\ze}{a_n}} =&
      \la{d(a)-1}\, \Da{d(a)-1}\,+\,\sum\limits_{\substack{c\in\sig0\\  c\neq a_{d(a)-1}}}\,\l c \, \D c \,+1\,,
\end{eqnarray}
and 
\begin{equation}\label{eq:recc:bound-n}
  \Da n\bydef\,\EE{\T{a_n}{a_n}}=0\,.
\end{equation}

Now, for $c\in\siga k, c\neq a_{k-1}$, the hitting time $\T{c}{a_n}$ of $a_n$ starting from $c$ is the sum of the time $t_c\bydef\T{c}{\p c}$ 
from $c$ to $\p c = a_k$ and $\T{a_k}{a_n}$, i.e. $\D c\,=\,t_c+\Da k$.\\ 
Using \reff{sommeprobasite} this leads to the following difference equation
\begin{equation}\label{eq:diff}
\left(\ma k\,+\,\la{k-1}\right)\,\Da k
\,=\,\ma k\,\Da{k+1}+\, \la{k-1}\,\Da{k-1}\,+\,\sum\limits_{\substack{c\in\siga k\\  c\neq a_{k-1}}}\,\l c \,t_c\,+\,1\,,
\end{equation}
with boundary conditions \reff{eq:recc:bound-n} and
\begin{eqnarray}
\m a\,\Da 0 \, =& \m a\,\Da1\,+\,\sum_{c\in\sig a}\,\l c \,t_c\,+\,1\,,\label{eq:recc:bound-a}\\
\la{d(a)-1}\,\Da{d(a)} =&
      \la{d(a)-1}\, \Da{d(a)-1}\,+\,\sum\limits_{\substack{c\in\sig0\\  c\neq a_{d(a)-1}}}\,\l c \, t_c \,+1\,,\label{eq:recc:bound-0}
\end{eqnarray}

\emph{Proof of \reff{T:enfant-parent}:} Suppose that  $0\leq j<n$.\\  
The idea is to write $\Da j$ as 
\begin{equation}
 \Da j\,=\,\Da0\,-\,\sum_{k=0}^{j-1}\,\left(\Da k \,-\, \Da{k+1} \right)\,.
\end{equation}
Together with boundary condition \reff{eq:recc:bound-n} which gives $\Da0\,=\,\sum_{k=0}^{n-1}\,\left(\Da k \,-\, \Da{k+1} \right)$, 
this leads to
\begin{equation}
 \Da j = \sum_{k=j}^{n-1}\,\left( \Da k-\,\Da{k+1} \right)\,.
\end{equation}
Let $\a k\bydef \dfrac{1}{\ma k}\,\Bigl( 1+\,\sum\limits_{\substack{c\in\siga k\\  c\neq a_{k-1}}}\,\l c \,t_c\Bigr)$ and 
$\b k \bydef \dfrac{\la{k-1}}{\ma k}$ for $k=1,\ldots,n-1$.\\
We can rewrite \reff{eq:diff} in the form 
\begin{equation}
 \Da k-\,\Da{k+1}\,=\,\b k\,(\Da {k-1}- \Da k)\,+\a k\,,
\end{equation}
and then iterate this difference equation down to $k=0$ (recall that $a_0=a$), to obtain
\begin{equation}
 \Da k-\,\Da{k+1}\,=\,\prod_{i=1}^k\,\b i\,(\Da 0- \Da 1)\,+\sum_{l=1}^k\,\a l\,\prod_{i=l+1}^k\,\b i\,.
\end{equation}
By condition \reff{eq:recc:bound-a}, we have
\begin{equation}
 \Da 0- \Da 1 \,=\, \dfrac{1}{\m a}\,\Bigl( 1+\,\sum_{c\in\sig a}\,\l c \,t_c\Bigr)\,\defby\,\a 0\,,
\end{equation}
Thus
\begin{equation}
 \Da k-\,\Da{k+1}\,=\,\sum_{l=0}^k\,\a l\,\prod_{i=l+1}^k\,\b i\,.
\end{equation}
And finally
\begin{equation}
 \Da j \,=\, \sum_{k=j}^{n-1}\,\left( \Da k-\,\Da{k+1} \right)\,=\, \sum_{k=j}^{n-1}\,\sum_{l=0}^k\,\a l\,\prod_{i=l+1}^k\,\b i \,.
\end{equation}
The quantities $\a l\,\prod_{i=l+1}^k\,\b i$ can be expressed in terms of the invariant measure $\pi$:\\
by definition \reff{invmes} of the i.p.m. $\pi$, we have
\begin{equation}\label{eq:produit-rapports}
 \dfrac{1}{\ma l}\,\prod_{i=l+1}^k\, \dfrac{\la{i-1}}{\ma i}
      \,=\, \dfrac{1}{\ma k}\,\prod_{i=l}^{k-1}\, \dfrac{\la i}{\ma i}
     \,=\, \dfrac{1}{\ma k}\,\dfrac{\pii{a_l}}{\pii{a_k}} \,,
\end{equation}
and using reversibility and expression \reff{ta} for $t_c\bydef \E{\T{c}{\p c}}$, we have
\begin{equation}
 \l c \,t_c\,=\,\dfrac{\l c\,\pii{\B c}}{\m c\,\pii c}\,=\,\dfrac{\pii{\B c}}{\pii{\p c}}\,=\,\dfrac{\pii{\B c}}{\pii{a_l}}\,,
\end{equation}
since by definition, $\p c=a_l$ for $c\in\siga l$. 
We get
\begin{eqnarray}
 \a 0\,\prod_{i=1}^k\,\b i 
    &=&\dfrac{1}{\m a}\,\prod_{i=1}^k\, \dfrac{\la{i-1}}{\ma i}\,\Bigl( 1+\,\sum_{c\in\sig a}\,\l c \,t_c\Bigr) \nonumber\\
    &=&\dfrac{1}{\ma k\,\pii{a_k}}\,\Bigl( \pii{a}\,+\,\sum_{c\in\sig a}\,\pii{\B c}\Bigr)\,=\,\dfrac{\pii{\B a}}{\ma k\,\pii{a_k}} \,,
\end{eqnarray}
and
\begin{eqnarray}
 \a l\,\prod_{i=l+1}^k\,\b i 
    &=& \dfrac{1}{\ma l}\,\prod_{i=l+1}^k\, \dfrac{\la{i-1}}{\ma i}\,
	  \Bigl( 1+\,\sum\limits_{\substack{c\in\siga l\\  c\neq a_{l-1}}}\,\l c \,t_c\Bigr) \nonumber\\
    &=& \dfrac{1}{\ma k\,\pii{a_k}}\,\Bigl( \pii{a_l}\,+\,\sum\limits_{\substack{c\in\siga l\\  c\neq a_{l-1}}}\,\,\pii{\B c}\Bigr) \,.
\end{eqnarray}
Thus
\begin{equation}
  \Da j \,=\, \sum_{k=j}^{n-1}\, \dfrac{1}{\ma k\,\pii{a_k}}\,
   \biggl(\, \sum_{l=1}^k\,\Bigl( \pii{a_l}\,+\,\sum\limits_{\substack{c\in\siga l\\  c\neq a_{l-1}}}\,\,\pii{\B c}\Bigr)\,+\,\pii{\B a} \,\biggr)\,,
\end{equation}
and equation \reff{T:enfant-parent} is proven, using the decomposition \reff{DecompBa} of $\Ba k$.\\


\emph{Proof of \reff{T:parent-enfant}:} Suppose that  $d(a)\geq j>n$.\\
The idea is now to iterate the difference equation in the other direction.\\ 
We write $\Da j$ as 
\begin{equation}
 \Da j\,=\,\Da{d(a)}\,-\,\sum_{k=j+1}^{d(a)}\,\left( \Da k-\,\Da{k-1} \right)\,,
\end{equation}
which with the boundary condition \reff{eq:recc:bound-n} gives
\begin{equation}
 \Da j\,=\,\,\sum_{k=n+1}^{j}\,\left( \Da k-\,\Da{k-1} \right)\,.
\end{equation}
We rewrite \reff{eq:diff} as
\begin{equation}
 \la{k-1}\,(\Da k- \Da {k-1})\,=\,\ma k\,(\Da{k+1}\,-\,\Da k)
	    \,+\,\Bigl( 1+\,\sum\limits_{\substack{c\in\siga k\\  c\neq a_{k-1}}}\,\l c \,t_c\Bigr)\,,
\end{equation}
which gives
\begin{equation}
 \Da k-\,\Da{k-1}\,=\,\bt k\,(\Da{k+1}- \Da  k)\,+\at k\,,
\end{equation}
with $\at k\bydef \dfrac{1}{\la{k-1}}\,\Bigl( 1+\,\sum\limits_{\substack{c\in\siga k\\  c\neq a_{k-1}}}\,\l c \,t_c\Bigr)$, and
$\bt k \bydef \dfrac{\ma k}{\la{k-1}}$ for $n<k<d(a)$.\\ 
Iteration up to $k=d(a)$ (recall $a_{d(a)}=0$), yields
\begin{equation}
  \Da k-\,\Da{k-1}\,=\,\prod_{i=k}^{d(a)-1}\,\bt i\,(\Da{d(a)}- \Da{d(a)-1})\,+\sum_{l=k}^{d(a)-1}\,\at l\,\prod_{i=k}^{l-1}\,\bt i\,.
\end{equation}
And  by the boundary condition \reff{eq:recc:bound-0}
\begin{equation}
 \Da{d(a)}- \Da{d(a)-1}\,=\,\dfrac{1}{\la{d(a)-1}}\,
		\Bigl( 1+\,\sum\limits_{\substack{c\in\siga{d(a)}\\  c\neq a_{d(a)-1}}}\,\l c \,t_c\Bigr)\,\defby\,\at{d(a)} \,.
\end{equation}
Thus
\begin{equation}
  \Da k-\,\Da{k-1}\,=\,\sum_{l=k}^{d(a)}\,\at l\,\prod_{i=k}^{l-1}\,\bt i\,,
\end{equation}
and then
\begin{equation}
  \Da j\,=\,\,\sum_{k=n+1}^{j}\,\,\sum_{l=k}^{d(a)}\,\at l\,\prod_{i=k}^{l-1}\,\bt i\,.
\end{equation}
As above this can be expressed in terms of the i.p.m. $\pi$, we get
\begin{eqnarray}
  \Da j&=&\sum_{k=n+1}^{j}\,\frac{1}{\la{k-1}\,\pii{a_k}}\,\sum_{l=k}^{d(a)}\,
	  \Bigl(\pii{a_l}\,+\,\sum\limits_{\substack{c\in\siga l\\  c\neq a_{l-1}}}\,\pii{\B c}\Bigr)\nonumber\\
      &=& \sum_{k=n+1}^{j}\,\frac{1-\pii{\Ba{k-1}}}{\la{k-1}\,\pii{a_k}} \\
      &=& \sum_{k=n}^{j-1}\,\frac{1-\pii{\Ba k}}{\ma k\,\pii{a_k}} \;.\nonumber
\end{eqnarray}
The second equality is due to the decomposition \reff{DecompBabar}, and the last one is obtained using reversibility.\\
Finally, equation \reff{T:parent-enfant} is proven by inverting $j$ and $n$ in this last formula.
\begin{flushright}
 $\Box$ 
\end{flushright}



\noindent{\bf{Proof of Proposition \ref{prop:moment2}}}\\
The proof follows the same lines as those of proposition \ref{prop:esp}. \\
We fix $a\in J$ and $0\leq n\leq d(a)$. The application of \reff{eq:claim} for $F(T)=T^2$ together with \reff{sommeprobasite} gives 
\begin{eqnarray}
  \EE{\T{a_k}{a_n}^2} &=& \ma k\, \EE{\T{a_{k+1}}{a_n}^2}\,+\,\la{k-1}\,\EE{\T{a_{k-1}}{a_n}^2}
		  \,+\,\sum\limits_{\substack{c\in\siga k\\  c\neq a_{k-1}}}\l c \,\EE{\T{c}{a_n}^2} \,+\,1\\
    && +\,2\,\biggl(\,\ma k\, \EE{\T{a_{k+1}}{a_n}}\,+\,\la{k-1}\,\EE{\T{a_{k-1}}{a_n}}
		  \,+\,\sum\limits_{\substack{c\in\siga k\\  c\neq a_{k-1}}}\l c \,\EE{\T{c}{a_n}}\,\biggr)\,,\nonumber
 \end{eqnarray}
and by \reff{eq:recc}, the second line is equal to $2\,\EE{\T{a_k}{a_n}}-2$. \\ 
We denote by $\De x$ the second moment of $\T{x}{a_n}$ and keep the notation $\D x$ for its expectation: 
\begin{equation}
 \De x\,=\,\EE{\T{x}{a_n}^2} \,,\,\text{ and} \quad \D x\,=\,\EE{\T{x}{a_n}}
\end{equation}
Thus for each $0<k<d(a)$, $k\neq n$, we have 
\begin{equation}\label{eq:recc2}
 \Dea k \,=\,\ma k\, \Dea{k+1}\,+\,\la{k-1}\,\Dea{k-1}
		  \,+\,\sum\limits_{\substack{c\in\siga k\\  c\neq a_{k-1}}}\l c \,\De c \,+\,2\,\Da k\,-\,1\,,
\end{equation}
and the same computations in sites $a_0=a$ and $a_{d(a)}=0$, gives the boundary conditions
\begin{eqnarray}
\Dea 0 \bydef& \EE{\T{a}{a_n}^2} =& \m a\,\Dea1\,+\,\sum_{c\in\sig a}\,\l c \,\De c\,+\,2\,\D a \,-\,1\,,\label{eq:recc2:bound-a}\\
\Dea{d(a)} \bydef& \EE{\T{\ze}{a_n}^2} =& \la{d(a)-1}\, \Dea{d(a)-1}
	\,+\,\sum\limits_{\substack{c\in\sig0\\  c\neq a_{d(a)-1}}}\,\l c \, \De c \,+\,2\,\D 0 \,-\,1\,,\label{eq:recc2:bound-0}
\end{eqnarray}
and for $k=n$
\begin{equation}\label{eq:recc2:bound-n}
  \Dea n\bydef\,\EE{\T{a_n}{a_n}^2}=0\,.
\end{equation}
As in the proof of Proposition \ref{prop:esp}, we can decompose $\T{c}{a_n}=\T{c}{\p c}+\T{a_k}{a_n}$, for $c\in\siga k$, $c\neq a_{k-1}$, 
and use \reff{sommeprobasite} to get the analogous difference equation of \reff{eq:diff} for second moments
\begin{eqnarray}\label{eq:diff2}
\left(\ma k\,+\,\la{k-1}\right)\,\Dea k &=& \ma k\,\Dea{k+1}+\, \la{k-1}\,\Dea {k-1}  \\
    	 && \,+\,\sum\limits_{\substack{c\in\siga k\\  c\neq a_{k-1}}}\,\l c \,\left(\,\EE{\T{c}{\p c}^2}\,+\,2\,t_c\,\Da k\,\right)
	  \,+\,2\,\Da k\,-\,1\,.\nonumber
\end{eqnarray}

\emph{Proof of \reff{T2:enfant-parent}:} When $0\leq j<n$, following the proof of \reff{T:enfant-parent} we get
\begin{equation}
 \Dea j \,=\,\sum_{k=j}^{n-1}\,\sum_{l=0}^k\,\a l\,\prod_{i=l+1}^k\,\b i \,,
\end{equation}
with $\b i \bydef \dfrac{\la{i-1}}{\ma i}$, $\a l\bydef \dfrac{1}{\ma l}\, \biggl( \;\sum\limits_{\substack{c\in\siga l\\  c\neq a_{l-1}}} \l c \,
\left(\,\EE{\T{c}{\p c}^2}\,+\,2\,t_c\,\Da l\,\right) \,+\,2\,\Da l-\,1\,\biggr)$ for $0<i,l<n$, 
and $\a 0\bydef \dfrac{1}{\m a}\, \biggl( \displaystyle\;\sum_{c\in\sig a}\,\l c \,\left(\,\EE{\T{c}{\p c}^2}
  \,+\,2\,t_c\,\D a\,\right) \,+\,2\,\D a\,-\,1 \,\biggr)$.

In order to get expression \reff{T2:enfant-parent}, we have now to express these quantities in terms of i.p.m. $\pi$. 
Let us start with the following lemma:
\begin{lem}\label{lem:alpha:prod:beta}
 \begin{equation}\label{eq:lem:alpha:l}
  \a l\,\prod_{i=l+1}^k\,\b i \,=\,\frac{1}{\ma k\,\pii{a_k}} \, \biggl(\, \sum\limits_{\substack{c\in\siga l\\  c\neq a_{l-1}}}\,
				\sum_{b\in\B c}\, \frac{2\,\pii{\B b}^2}{\m b\,\pii{b}} \,+\, \Bigl(\,2\,\Da l \,-\,1 \,\Bigr) 
			      \, \Bigl(\,\pii{\Ba l}\,-\,\pii{\Ba{l-1}} \,\Bigr) \,\biggr) 
 \end{equation}
 \begin{equation}\label{eq:lem:alpha:0}
    \a 0\,\prod_{i=1}^k\,\b i \,=\\\frac{1}{\ma k\,\pii{a_k}}\, \biggl(\, \sum_{c\in\sig a}\,\sum_{b\in\B c}\, \frac{2\,\pii{\B b}^2}{\m b\,\pii{b}}
      \,+\, \bigl(\,2\,\D a \,-\,1 \,\bigr) \, \pii{\B a} \,\biggr)
 \end{equation}
\end{lem}

\noindent{\bf{Proof of Lemma \ref{lem:alpha:prod:beta}}}\\ 
Let $0<l<n$. We recall that by \reff{eq:produit-rapports}, we have
\begin{equation}\label{eq:produit-rapports:bis}
 \dfrac{1}{\ma l}\,\prod_{i=l+1}^k\, \dfrac{\la{i-1}}{\ma i} \,=\,\dfrac{\pii{a_l}}{\ma k\,\pii{a_k}} \,.
\end{equation}
It remains to compute the term in the parenthesis in the definition of $\a l$.\\
Using reversibility, $\dfrac{\l c}{\m c\,\pii{c}}=\dfrac{1}{\pii{a_l}}$ for $c\in\siga l$, and then 
by equation \reff{eq:lem:ub}. we have
\begin{eqnarray}\label{eq:proof:lem:alpha:l}
&&\sum\limits_{\substack{c\in\siga l\\  c\neq a_{l-1}}}\,\l c \, \left(\,\EE{\T{c}{\p c}^2}\,+\,2\,t_c\,\Da l\,\right)
	 \,+\,2\,\Da l\,-\,1  \nonumber\\
=&& \frac{1}{\pii{a_l}} \,\sum\limits_{\substack{c\in\siga l\\  c\neq a_{l-1}}}\,\sum_{b\in\B c}\, \frac{2\,\pii{\B b}^2}{\m b\,\pii{b}}\,+\,
   \bigl(\,2\,\Da l \,-\,1 \,\bigr) \, \bigl(\,1\,+\,\sum\limits_{\substack{c\in\siga l\\ c\neq a_{l-1}}}\,\l c \,t_c \,\bigr)\\
=&& \frac{1}{\pii{a_l}} \, \biggl(\, \sum\limits_{\substack{c\in\siga l\\  c\neq a_{l-1}}}\,\sum_{b\in\B c}\, \frac{2\,\pii{\B b}^2}{\m b\,\pii{b}}
      \,+\, \bigl(\,2\,\Da l \,-\,1 \,\bigr) \, \Bigl(\,\pii{\Ba l}\,-\,\pii{\Ba{l-1}} \,\Bigr) \,\biggr) \,,\nonumber
\end{eqnarray}
since by \reff{ta}
\begin{eqnarray}
 \bigl(\,1\,+\,\sum\limits_{\substack{c\in\siga l\\  c\neq a_{l-1}}}\,\l c \,t_c  \,\bigr) 
	&=& \frac{1}{\pii{a_l}}\,\Bigl(\,\pii{a_l}\,+\,\sum\limits_{\substack{c\in\siga l\\  c\neq a_{l-1}}}\,\pii{\B c}  \,\Bigr)  \\
	&=& \frac{1}{\pii{a_l}}\,\Bigl(\,\pii{\Ba l}\,-\,\pii{\Ba{l-1}} \,\Bigr) \,. \nonumber
\end{eqnarray}
Equation \reff{eq:lem:alpha:l} is the consequence of  \reff{eq:produit-rapports:bis} and \reff{eq:proof:lem:alpha:l}.\\
For the proof of \reff{eq:lem:alpha:0} we use similar arguments for $l=0$:
\begin{equation}
 \dfrac{1}{\m a}\,\prod_{i=1}^k\, \dfrac{\la{i-1}}{\ma i} \,=\,\dfrac{\pii{a}}{\ma k\,\pii{a_k}} \,,
\end{equation}
and
\begin{eqnarray}
&& \sum_{c\in\sig a}\,\l c \,\left(\,\EE{\T{c}{\p c}^2} \,+\,2\,t_c\,\D a\,\right) \,+\,2\,\D a\,-\,1 \nonumber\\
=&& \frac{1}{\pii{a}} \,\sum_{c\in\sig a}\,\sum_{b\in\B c}\, \frac{2\,\pii{\B b}^2}{\m b\,\pii{b}}\,+\,
   \bigl(\,2\,\D a \,-\,1 \,\bigr) \, \bigl(\,1\,+\,\sum_{c\in\sig a}\,\l c \,t_c \,\bigr)\\
=&& \frac{1}{\pii{a}} \, \biggl(\, \sum_{c\in\sig a}\,\sum_{b\in\B c}\, \frac{2\,\pii{\B b}^2}{\m b\,\pii{b}}
      \,+\, \bigl(\,2\,\D a \,-\,1 \,\bigr) \, \pii{\B a} \,\biggr) \,.\nonumber
\end{eqnarray}
\begin{flushright}
 $\Box$ 
\end{flushright}

We now have 
\begin{equation}
 \Dea j \,=\,\sum_{k=j}^{n-1}\,\frac{1}{\ma k\,\pii{a_k}}\,\bigl(\, \sum_{l=1}^k\,\Gamma_l \,+\, \Gamma_0 \,\bigr) \,,
\end{equation}
with
\begin{equation}\label{eq:gamma l}
 \Gamma_l \bydef \,\sum\limits_{\substack{c\in\siga l\\  c\neq a_{l-1}}}\,
				\sum_{b\in\B c}\, \frac{2\,\pii{\B b}^2}{\m b\,\pii{b}} \,+\, \Bigl(\,2\,\Da l \,-\,1 \,\Bigr) 
			      \, \Bigl(\,\pii{\Ba l}\,-\,\pii{\Ba{l-1}} \,\Bigr)
\end{equation}
and 
\begin{equation}\label{eq:gamma 0}
 \Gamma_0 \bydef \,\sum_{c\in\sig a}\,\sum_{b\in\B c}\, \frac{2\,\pii{\B b}^2}{\m b\,\pii{b}}
      \,+\, \bigl(\,2\,\D a \,-\,1 \,\bigr) \, \pii{\B a} \,.
\end{equation}
Summing over $l$, the second term in \reff{eq:gamma l}, we get
\begin{eqnarray}\label{eq:sumBl-Bl-1}
&& \sum_{l=1}^k\, \Bigl(\,2\,\Da l \,-\,1 \,\Bigr) \, \Bigl(\,\pii{\Ba l}\,-\,\pii{\Ba{l-1}} \,\Bigr) \nonumber\\
&=& \sum_{l=1}^k\, \Bigl(\,2\,\Da l \,-\,1 \,\Bigr) \,\pii{\Ba l} 
	      \,-\,\sum_{l=1}^k\, \Bigl(\,2\,\Da l \,-\,1 \,\Bigr) \,\pii{\Ba{l-1}} \nonumber\\
&=& \sum_{l=1}^k\, \Bigl(\,2\,\bigl(t_{a_l}+\Da{l+1} \bigr)\,-\,1 \,\Bigr) \,\pii{\Ba l} 
	      \,-\,\sum_{l=0}^{k-1}\, \Bigl(\,2\,\Da{l+1} \,-\,1 \,\Bigr) \,\pii{\Ba l} \\
&=& \sum_{l=1}^k\,2\,t_{a_l}\,\pii{\Ba l}  \,+\,\sum_{l=1}^k\, \Bigl(\,2\,\Da{l+1} \,-\,1 \,\Bigr) \,\pii{\Ba l} 
	      \,-\,\sum_{l=0}^{k-1}\, \Bigl(\,2\,\Da{l+1} \,-\,1 \,\Bigr) \,\pii{\Ba l}\nonumber\\
&=& \sum_{l=1}^k\,\frac{2\,\pii{\Ba l}^2}{\ma l\,\pii{a_l}}  \,+\, \bigl(\, 2\,\Da{k+1} -\,1\,\bigr)\, \pii{\Ba k} 
	      \,-\,\bigl(\,2\,\Da 1 \,-\,1 \,\bigr) \,\pii{\B a} \,,\nonumber
\end{eqnarray}
where we have used that for $l<n$, we have $\T{a_l}{a_n}\,=\,\T{a_l}{a_{l+1}}\,+\,\T{a_{l+1}}{a_n}\,$.\\
And the boundary condition \reff{eq:recc:bound-a} gives
\begin{eqnarray}
&& \pii{\B a}\, \Bigl(\bigl(\,2\,\D a \,-\,1 \,\bigr)\,-\,\bigl(\,2\,\Da 1 \,-\,1 \,\bigr) \,\Bigr)  \nonumber\\
&=& 2\,\pii{\B a}\,\bigl(\,\D a \,-\,\Da 1 \,\bigr) \nonumber\\
&=& \frac{2\,\pii{\B a}}{\m a}\, \bigl(\,1\,+\,\sum_{c\in\sig a}\,\l c \,t_c \,\bigr) \nonumber\\
&=& \frac{2\,\pii{\B a}}{\m a\,\pii{a}}\, \bigl(\,\pii{a}\,+\,\sum_{c\in\sig a}\,\pii{\B c} \,\bigr) \nonumber\\
&=& \frac{2\,\pii{\B a}^2}{\m a\,\pii{a}}\,. \nonumber
\end{eqnarray}
Thus
\begin{equation}
 \Gamma_0 \,-\,\bigl(\,2\,\Da 1 \,-\,1 \,\bigr) \,\pii{\B a}\,=\,\sum_{c\in\sig a}\,\sum_{b\in\B c}\, \frac{2\,\pii{\B b}^2}{\m b\,\pii{b}}
      \,+\,\frac{2\,\pii{\B a}^2}{\m a\,\pii{a}} \,=\,\sum_{b\in\B a}\, \frac{2\,\pii{\B b}^2}{\m b\,\pii{b}} \,,
\end{equation}
since by definition $\displaystyle \B a=\{a\}\bigcup_{c\in\sig a}\,\B c$ .\\
Hence
\begin{eqnarray}
\sum_{l=1}^k\,\Gamma_l \,+\, \Gamma_0 &=& \sum_{l=1}^k\,\Biggl(\,\sum\limits_{\substack{c\in\siga l\\  c\neq a_{l-1}}}\,
	\sum_{b\in\B c}\, \frac{2\,\pii{\B b}^2}{\m b\,\pii{b}} \,+\,\frac{2\,\pii{\Ba l}^2}{\ma l\,\pii{a_l}} \Biggr)
 \,+\,\sum_{b\in\B a}\, \frac{2\,\pii{\B b}^2}{\m b\,\pii{b}} \\
&&\quad + \, \bigl(\, 2\,\Da{k+1} -\,1\,\bigr)\, \pii{\Ba k} \\
&=&   \sum_{b\in\Ba k}\, \frac{2\,\pii{\B b}^2}{\m b\,\pii{b}} \, + \, \bigl(\, 2\,\Da{k+1} -\,1\,\bigr)\, \pii{\Ba k} \,,
\end{eqnarray}
where the second equality is due to the decomposition \reff{DecompBa} of $\Ba k$.\\
Finally, 
\begin{equation}
 \Dea j \,=\,
\sum_{k=j}^{n-1}\,\frac{1}{\ma k\,\pii{a_k}}\,\Biggl(\, \sum_{b\in\Ba k}\, \frac{2\,\pii{\B b}^2}{\m b\,\pii{b}} 
	  \, + \,2\,\pii{\Ba k}\,\Da{k+1} \,-\, \pii{\Ba k} \,\Biggr)    \,,
\end{equation}
and equation \reff{T2:enfant-parent} is obtained using \reff{T:enfant-parent}.\\

\emph{Proof of \reff{T2:parent-enfant}:} Suppose that  $d(a)\geq j>n$.\\
As in the proof of \reff{T:parent-enfant}, from the difference equation \reff{eq:diff2}, we obtain
\begin{equation}
  \Dea j\,=\,\,\sum_{k=n+1}^{j}\,\,\sum_{l=k}^{d(a)}\,\at l\,\prod_{i=k}^{l-1}\,\bt i\,,
\end{equation}
with $\at l\bydef \dfrac{1}{\la{l-1}}\, \biggl( \;\sum\limits_{\substack{c\in\siga l\\  c\neq a_{l-1}}} \l c \,
\left(\,\EE{\T{c}{\p c}^2}\,+\,2\,t_c\,\Da l\,\right) \,+\,2\,\Da l-\,1\,\biggr)$ for $n<l\leq d(a)$,\\
and $\bt i \bydef \dfrac{\ma i}{\la{i-1}}$ for $n<i<d(a)$.\\
By the same arguments as those of the proof of Lemma \ref{lem:alpha:prod:beta}, 
we get $\displaystyle \at l\,\prod_{i=k}^{l-1}\,\bt i \,=\, \frac{1}{\la{k-1}\,\pii{a_k}} \; \Gamma_l$, 
where $\Gamma_l$ is defined by equation \reff{eq:gamma l}. And thus, using reversibility
\begin{equation}
  \Dea j\,=\,\,\sum_{k=n+1}^{j}\,\frac{1}{\la{k-1}\,\pii{a_k}}\;\sum_{l=k}^{d(a)}\;\Gamma_l
	  \,=\,\,\sum_{k=n}^{j-1}\,\frac{1}{\ma k\,\pii{a_k}}\;\sum_{l=k+1}^{d(a)}\;\Gamma_l\,.
\end{equation}
It remains to compute $\displaystyle\sum_{l=k+1}^{d(a)}\;\Gamma_l$. First, by definition \reff{DecompCa}, we have
\begin{equation}
 \sum_{l=k+1}^{d(a)}\,\sum\limits_{\substack{c\in\siga l\\  c\neq a_{l-1}}}\,\sum_{b\in\B c}\, \frac{2\,\pii{\B b}^2}{\m b\,\pii{b}}
	\,=\, \sum_{c\in\Ca k}\, \frac{2\,\pii{\B c}^2}{\m c\,\pii{c}} \,.
\end{equation}
For the sum of the second term in \reff{eq:gamma l}, we use a computation similar to \reff{eq:sumBl-Bl-1} with the following decomposition: 
for $l>n$, we have $\T{a_{l+1}}{a_n}\,=\,\T{a_{l+1}}{a_l}\,+\,\T{a_l}{a_n}\,$, and thus
\begin{eqnarray}
&& \sum_{l=k+1}^{d(a)}\, \Bigl(\,2\,\Da l \,-\,1 \,\Bigr) \, \Bigl(\,\pii{\Ba l}\,-\,\pii{\Ba{l-1}} \,\Bigr) \nonumber\\
&=& \sum_{l=k+1}^{d(a)}\, \Bigl(\,2\,\Da l \,-\,1 \,\Bigr) \,\pii{\Ba l} 
	 \,-\,\sum_{l=k}^{d(a)-1}\, \Bigl(\,2\,\bigl(\EE{\T{a_{l+1}}{a_l}} + \Da{l} \bigr) \,-\,1 \,\Bigr) \,\pii{\Ba{l}} \\
&=& \Bigl(\, 2\,\Da{d(a)} -\,1\,\Bigr)
	 \,-\,\Bigl(\,2\,\Da{k} \,-\,1 \,\Bigr)\,\pii{\Ba{k}} \,-\,2\, \sum_{l=k}^{d(a)-1}\,\EE{\T{a_{l+1}}{a_l}} \,\pii{\Ba{l}} \,,\nonumber
\end{eqnarray}
since by definition $\Ba{d(a)}=\B0=I$. \\
We now write $\Da{d(a)} \,\bydef\,\EE{\T{a_{d(a)}}{a_n}}$ as $\EE{\T{a_{d(a)}}{a_k}} \,+\, \Da k$, and since by \reff{T:parent-enfant}, we have
\begin{equation}
 \EE{\T{a_{d(a)}}{a_k}} \,=\,\sum_{l=k}^{d(a)-1}\,\EE{\T{a_{l+1}}{a_l}} \,=\,\sum_{l=k}^{d(a)-1}\,\frac{1-\pii{\Ba l}}{\ma l\,\pii{a_l}} \,, 
\end{equation}
we get
\begin{eqnarray}
&& \sum_{l=k+1}^{d(a)}\, \Bigl(\,2\,\Da l \,-\,1 \,\Bigr) \, \Bigl(\,\pii{\Ba l}\,-\,\pii{\Ba{l-1}} \,\Bigr) \nonumber\\
&=& \Bigl(\,2\,\Da{k} \,-\,1 \,\Bigr)\,\Bigl(\,1\,-\,\pii{\Ba{k}}\,\Bigr) 
	  \,+\,2\, \sum_{l=k}^{d(a)-1}\,\frac{\Bigl(\,1\,-\,\pii{\Ba l}\,\Bigr)^2}{\ma l\,\pii{a_l}} \\
&=& \Bigl(\,2\,\Da{k} \,-\,1 \,\Bigr)\,\pii{\,\Babar{k}}
	  \,+\,2\, \sum_{l=k}^{d(a)-1}\,\frac{\pii{\,\Babar l}^2}{\ma l\,\pii{a_l}} \,.\nonumber
\end{eqnarray}
Finally
\begin{equation}
 \Dea j =\,\sum_{k=n}^{j-1}\,\frac{1}{\ma k\,\pii{a_k}}\,
		  \Biggl(\,\sum_{c\in\Ca{k}}\, \frac{2\,\pii{\B c}^2}{\m c\,\pii{c}} 
			  \,+\,2\, \sum_{l=k}^{d(a)-1}\,\frac{\pii{\,\Babar l}^2}{\ma l\,\pii{a_l}} 
			   \,+\,\Bigl(\,2\,\Da{k} \,-\,1 \,\Bigr)\,\pii{\,\Babar{k}} \,\Biggr) \,,
\end{equation}
and equation \reff{T2:parent-enfant} is obtained by inverting $j$ and $n$, using expression \reff{T:parent-enfant} for $\EE{\T{a_n}{a_j}}$, 
and since by definition \reff{def:la}, the path $\ell(a_k)$ is equal to $\cup_{l=k}^{d(a)-1}\{a_l\}$.

\begin{flushright}
 $\Box$ 
\end{flushright}


\begin{thebibliography}{10}


\bibitem{AldousRW}
Aldous, D.: Random walks on finite groups and rapidly mixing {M}arkov chains.
\newblock In: S\'eminaire de probabilit\'es XVII, \emph{Lecture Notes in Math.} vol.
  986, pp. 243--297. Springer, Berlin (1983)

\bibitem{DiaAlShCST}
Aldous, D., Diaconis, P.: Shuffling cards and stopping times.
\newblock Amer. Math. Monthly \textbf{93}(5), 333--348 (1986)

\bibitem{DiaAlSUT}
Aldous, D., Diaconis, P.: Strong uniform times and finite random walks.
\newblock Adv. in Appl. Math. \textbf{8}(1), 69--97 (1987)

\bibitem{AldousFill}
Aldous, D., Fill J.: Reversible Markov Chains and Random Walks on Graphs.
\newblock Monograph in preparation \emph{http://www.stat.berkeley.edu/users/aldous/RWG/book.html} 

\bibitem{barberfer09}
Barrera, J., Bertoncini, O., Fern\'andez, R.: Abrupt Convergence and Escape Behavior for Birth and Death Chains
\newblock J. Stat. Phys. \textbf{137}, 595-623 (2009)

\bibitem{olivphd}
Bertoncini, O.: Convergence abrupte et m\'etastabilit\'e.
\newblock Ph.D. thesis, Universit\'e de Rouen (2007)

\bibitem{berbarfer08}
Bertoncini, O., Barrera, J., Fern\'andez, R.: Cut-off and exit from metastability:
  two sides of the same coin.
\newblock C. R. Acad. Sci. Paris, Ser. I \textbf{346}, 691--696 (2008)

\bibitem{CGOV}
Cassandro, M., Galves, A., Olivieri, E., Vares, M.E.: Metastable behavior of
  stochastic dynamics: a pathwise approach.
\newblock J. Statist. Phys. \textbf{35}(5-6), 603--634 (1984)


\bibitem{DorGovMend08}
Dorogovtsev, S.N., Goltsev, A. V., Mendes, J.F.F.: 
\newblock Critical phenomena in complex networks.
\newblock Rev. Mod. Phys., 80, 1275 (2008)

\bibitem{DoroLectures10}
Dorogovtsev, S.N.:
\newblock Lectures on Complex Networks.
\newblock Oxford University Press, Oxford (2010)

\bibitem{Lyons90}
Lyons, R.: Random Walks and Percolation on Trees.
\newblock Ann. Probab. \textbf{18}(3), 931--958 (1990)

\bibitem{Woess}
Woess, W.: Denumerable {M}arkov Chains. Generating functions, boundary theory, random walks on trees.
\newblock EMS Textbooks in Mathematics. European Mathematical Society (EMS), Z\"{u}rich (2009)

\bibitem{YcartMart}
Mart{\'{\i}}nez, S., Ycart, B.: Decay rates and cutoff for convergence and
  hitting times of {M}arkov chains with countably infinite state space.
\newblock Adv. in Appl. Probab. \textbf{33}(1), 188--205 (2001)


\end{thebibliography}
\end{document}